\documentstyle[amssymb,12pt]{amsart}

\newtheorem{lem}{Lemma}
\newtheorem{prop}{Proposition}

\theoremstyle{definition}

\theoremstyle{definition}

\theoremstyle{remark}
\newtheorem{rem}{Remark}

\numberwithin{equation}{section}

\begin{document}

\newcommand{\thmref}[1]{Theorem~\ref{#1}}
\newcommand{\secref}[1]{Sect.~\ref{#1}}
\newcommand{\lemref}[1]{Lemma~\ref{#1}}
\newcommand{\propref}[1]{Proposition~\ref{#1}}
\newcommand{\corref}[1]{Corollary~\ref{#1}}
\newcommand{\remref}[1]{Remark~\ref{#1}}
\newcommand{\nc}{\newcommand}
\nc{\on}{\operatorname}
\nc{\ch}{\mbox{ch}}
\nc{\Z}{{\Bbb Z}}
\nc{\C}{{\Bbb C}}
\nc{\pone}{{\Bbb C}{\Bbb P}^1}
\nc{\pa}{\partial}
\nc{\F}{{\cal F}}
\nc{\arr}{\rightarrow}
\nc{\larr}{\longrightarrow}
\nc{\al}{\alpha}
\nc{\ri}{\rangle}
\nc{\lef}{\langle}
\nc{\W}{{\cal W}}
\nc{\la}{\lambda}
\nc{\ep}{\epsilon}
\nc{\su}{\widehat{\frak{sl}}_2}
\nc{\sw}{\frak{sl}}
\nc{\g}{\frak{g}}
\nc{\h}{\frak{h}}
\nc{\n}{\frak{n}}
\nc{\N}{\widehat{\n}}
\nc{\ab}{\frak{a}}
\nc{\G}{\widehat{\g}}
\nc{\De}{\Delta_+}
\nc{\gt}{\widetilde{\g}}
\nc{\Ga}{\Gamma}
\nc{\one}{{\bold 1}}
\nc{\hh}{\widehat{\h}}
\nc{\z}{{\frak Z}}
\nc{\zz}{{\cal Z}}
\nc{\Hh}{{\cal H}}
\nc{\qp}{q^{\frac{k}{2}}}
\nc{\qm}{q^{-\frac{k}{2}}}
\nc{\La}{\Lambda}
\nc{\wt}{\widetilde}
\nc{\qn}{\frac{[m]_q^2}{[2m]_q}}
\nc{\cri}{_{\on{cr}}}
\nc{\k}{h^\vee}
\nc{\sun}{\widehat{\sw}_N}
\nc{\HH}{{\cal H}_h(\sw_N)}
\nc{\ca}{\wt{{\cal A}}_{h,k}(\sw_2)}
\nc{\si}{\sigma}
\nc{\gl}{\widehat{\frak{g}\frak{l}}_2}
\nc{\el}{\ell}
\nc{\s}{s}
\nc{\Res}{\on{Res}}
\nc{\Vir}{{\cal Vir}}
\nc{\wz}{\left( \frac{w}{z} \right)}
\nc{\ST}{S_{TT}}
\nc{\Rep}{\on{Rep}}
\nc{\dzz}{\frac{dz}{z}}
\nc{\dww}{\frac{dw}{w}}
\nc{\bi}{\bibitem}

\title[Towards Deformed Chiral Algebras]{Towards Deformed Chiral
Algebras}\thanks{Contribution to the Proceedings of the Quantum Group
Symposium at the XXI$^{\on{th}}$ International Colloquium on Group
Theoretical Methods in Physics, Goslar 1996}

\author{Edward Frenkel}

\address{Department of Mathematics, Harvard University, Cambridge, MA
02138, USA}

\author{Nikolai Reshetikhin}

\address{Department of Mathematics, University of California, Berkeley, CA
94720, USA}

\maketitle

\section{Introduction}

\subsection{}
In \cite{LP1,FR,A1,FF:w,A2,L1} certain deformations of the Virasoro and
$\W$--algebras have been defined together with their free field
realizations. While the ordinary $\W$--algebras are symmetries of conformal
field theories \cite{BPZ,BS}, the deformed $\W$--algebras appear to be the
dynamical symmetries of the RSOS integrable models \cite{LP2,M1}. But the
deformed $\W$--algebras are still rather poorly understood. One can define
them in terms of generators and relations, but the relations are unwieldy
and not very illuminating. This makes the study of representation theory of
deformed $\W$--algebras very complicated (even in the case of the deformed
Virasoro algebra). One can hope that there is a new algebraic structure
behind the deformed $\W$--algebras, which would make them more
tractable. In this note we attempt to describe such a structure.

The current situation with deformed $\W$--algebras is actually similar
to the situation with the ordinary $\W$--algebras. If one tries to
describe them in terms of the ``Fourier coefficients of currents'',
one also runs into unwieldy formulas (except in the case of the
Virasoro algebra, which happens to be a Lie algebra). The right idea
is to consider instead of the Fourier coefficients of the generating
fields $W_i(z)$, the fields themselves. These fields together with
their derivatives and normally ordered products form a vector space
$V$ with the operation of operator product expansion (OPE),
i.e., whenever $A(z)$ and $B(z)$ belong to $V$, one can expand their
composition as
$$A(z) B(w) = \sum_{n \in \Z} C_n(w) (z-w)^n,$$ where all the
coefficients $C_n(w)$ belong to $V$. Thus we obtain a vertex operator
algebra (VOA) -- we recall the precise definition of VOA in
Sect.~2. This VOA can be described explicitly as a subalgebra of a
simple Heisenberg VOA \cite{FF:laws} or as the result of the quantum
Drinfeld-Sokolov reduction \cite{FF:ds}. The commutation relations of
the Fourier coefficients of the fields can in principle be computed
from the OPEs, but this is not necessary as the VOA structure suffices
for most purposes.

\subsection{}
In Sect.~3 we apply the same idea in the deformed case and define what we
call the deformed chiral algebra (the term ``chiral algebra'' is synonymous
to ``vertex operator algebra'', but we prefer the terminology ``deformed
chiral algebra'' to ``deformed VOA''). This definition has several
important differences with the undeformed case.

First of all, in the case of a VOA $V$, one assumes from the beginning that
each field $A(z)$ from $V$ acts on $V$, i.e., each Fourier coefficient of
$A(z)$ can be considered as a linear operator on $V$. In other words, each
VOA is manifestly a representation of itself. There are many indications
that this should not be so in the deformed case. Therefore, we distinguish
from the beginning the space of fields $V$ and the space of states $W$,
on which the fields act.

Second, while in the undeformed case the analytic continuation of the
composition of two fields $A(z)$ and $B(w)$ has poles only on the diagonal
$z=w$, in the deformed case we allow poles on various shifted diagonals
$z=w\gamma$, where $\gamma$ lies in a certain lattice in $\C^\times$. We
therefore consider the OPEs at each of these poles.

Finally, the locality axiom, which in the undeformed case states that
analytic continuations of $A(z) B(w)$ and $B(w) A(z)$ coincide, has to be
modified. Instead, we require that they differ by an action of an operator
$S(w/z)$ on $V \otimes V$ (the $S$--matrix). We call this property
$S$--commutativity. Note that a simple example of $S$--commutativity
appears already in super-VOAs, where $S$ acts as $\pm 1$ depending on the
parity of $A$ and $B$.

Thus, the main ingredients of a deformed chiral algebra (DCA) are operator
product expansions at shifted diagonals, and $S$--commutativity.

In Sect.~4 we give some examples of DCAs. One is associated to a Heisenberg
algebra, and the other is its subalgebra associated with the deformed
Virasoro algebra. Analogously, one can construct the DCAs corresponding to
other $\W$--algebras and to the affine quantum algebras $U_p(\sun)$ using
their free field realizations \cite{FF:w,A2,FJ,AOS}. The free field
realizations remain at the moment the only way known to us for producing
explicit examples of deformed chiral algebras.

We should say from the outset that we view our definition of DCA as just a
first step towards understanding the true meaning of the deformed
$\W$--algebras. One drawback of our definition is that in order to satisfy
the axioms we have to consider a very large space of fields. It seems that
many of them are redundant, but so far we have not found a natural way to
truncate this space.

\subsection{}
One of our motivations for introducing the DCAs associated to the deformed
$\W$--algebras is to understand better their connection with quantum affine
algebras (see Sect.~5).

Recall that the classical $\W$--algebra $\W(^L\g)$, where $^L\g$ is the
Langlands dual of $\g$, is isomorphic to the center $Z(\g)$ of
$U(\G)_{\on{cr}}$ -- the completed quantized universal enveloping algebra
of $\G$ at the critical level \cite{FF:ds}. A deformed classical
$\W$--algebras was defined in \cite{FR} as the center $Z_p(\g)$ of
$U_p(\G)_{\on{cr}}$ -- the completed quantized universal enveloping algebra
of $\G$ at the critical level. The deformed quantum $\W$--algebra
$Z_{p,q}(\g)$ is the quantization of the center considered as a Poisson
algebra. The quantization parameter is $1-p/q$, i.e., $Z_p(\g)$ is the
classical limit of $Z_{p,q}(\g)$ as $q \arr p$.

A construction of \cite{RS} associates to each finite-dimensional
representation $U$ of $U_p(\G)$ (i.e., of level $0$) a generating series
$t_U(z)$ of central elements from $Z_p(\g)$. In particular, the $i$th
fundamental representation of $U_p(\sun)$ corresponds to a series $t_i(z)$,
whose quantization gives rise to one of the generating fields of
$Z_{p,q}(\sw_N) \equiv \W_{p,q}(\sw_N)$ \cite{FF:w,A2}. The question is
what are the quantizations of central elements corresponding to other
(higher spin) irreducible representations of $U_p(\sun)$.

The DCA structure helps to answer this question. In the case of the
deformed Virasoro algebra, the generating field $T(w)$ corresponds to the
two-dimensional representation of $U_p(\su)$. We show in Sect.~6 that the
three--dimensional representation corresponds to the residue in the OPE of
$T(z)$ and $T(w)$ at $z=wq$. More generally, the field $T_n(w)$
corresponding to the $(n+1)$--dimensional representation of $U_p(\su)$ can
be defined inductively as the residue of the OPE of $T(z)$ and $T_{n-1}(w)$
at $z=wq^{n-1}$. At the end of Sect.~6 we derive various formal power
series relations between the Fourier coefficients of the fields $T_n(z)$.

The research of E.~Frenkel was supported by grants from the Packard
Foundation, the Sloan Foundations, and the NSF. The research of
N.~Reshetikhin was supported by an NSF grant.

\section{Vertex Operator Algebras}

In this section we recall the axioms of vertex operator algebra following
\cite{FKRW} (see also \cite{Kac} and \cite{FF:laws}, Sect. 4.1). These
axioms are equivalent to the original axioms \cite{Bor,FLM}; for the proof
of this equivalence, see \cite{Kac}.

Let $V=\oplus_{n=0}^{\infty}V_n\/$ be a $\Bbb Z_{+}\/$-graded vector space,
where $\dim V_n < \infty\/$ for all $n\/$. A field on $V\/$ of conformal
dimension $\Delta \in \Bbb Z\/$ is a power series $\phi (z) = \sum_{j \in
\Bbb Z} \phi_j z^{-j - \Delta}\/$, where $\phi_j \in \mbox{End}\, V \/$ and
$\phi_j V_n \subset V_{n-j}\/$.

If $z \in \C^\times$ is a non-zero complex number, then $\phi(z)$ can be
considered as a linear operator $V \arr \bar{V}$, where
$\bar{V}=\prod_{n=0}^\infty V_n$.

We have a natural pairing $\langle,\rangle: V_n^{*} \times V_n
\rightarrow \C$. A linear operator $P: V \rightarrow
\bar{V}$ can be represented by a set of finite-dimensional linear
operators $P^{j}_{i}: V_i \arr V_j, i,j \in
\Z$, such that $\langle A, P \cdot
B\rangle = \langle A, P^{j}_{i} \cdot B \rangle$ for $B \in V_i, A
\in V_j^{*}$. Let $P, Q$ be two linear operators $V \rightarrow
\bar{V}$. We say that the composition $PQ$ exists, if for any
$i,k \in \Z, B \in V_i, A \in V_k^*$ the series $\sum_{j \in \Z}
\langle A, P^{k}_{j} Q^{j}_{i} \cdot B \rangle$ converges absolutely.

Two fields $\phi (z)\/$ and $\psi(z)\/$
are called $\on{local}$ with respect to each other, if

\begin{itemize}
\item for any $z, w \in \C^\times$, such
that $|z|>|w|$, the composition $\phi(z) \psi(w)$ exists and can be
analytically continued to a rational operator-valued function on
$(\C^{\times})^{2}\backslash diagonal$, $R(\phi(z) \psi(w))$;

\item for any $z, w \in \C^\times$, such that $|w|>|z|$, the
composition $\psi(w) \phi(z)$ exists and can be analytically continued to a
rational operator-valued function on $(\C^{\times})^{2}\backslash
diagonal$, $R(\psi(w) \phi(z))$;

\item $R(\phi(z) \psi(w)) = R(\psi(w) \phi(z))$.
\end{itemize}

\vspace*{4mm}
\noindent {\bf Definition.}  A VOA algebra is a collection of the following
data: 

\begin{itemize}
\item a $\Z_+$--graded vector space $V$.

\item a linear map $Y : V \longrightarrow \on{End}\, V
\left[\left[z,z^{-1}\right]\right]\/$ which associates to each $A \in
V_n\/$ a field of conformal dimension $n\/$ $$Y(A,z) = \sum_{j \in \Bbb
Z}A_j z^{-j-n}\/.$$

\item A linear operator $T: V \arr V$ of degree $1$.

\item A vector $\Omega \in V$.

\end{itemize}

These data should satisfy the following axioms:

\begin{itemize}

\item[(A1)]
All fields $Y(A,z)\/$ are local with respect to each other.

\item[(A2)]
$T \Omega = 0$ and $\partial_z Y(A,z) = [T,Y(A,z)]$

\item[(A3)]
$Y(\Omega,z) = \on{Id}$ and $\lim_{z
\rightarrow 0} Y(A,z)|0\rangle = A\/$ for all $A \in V$.

\end{itemize}

\begin{prop}    \label{assoc}
A VOA $V$ automatically satisfies the following property: for any $A, B \in
V$,
\begin{equation}
R(Y(A,z) Y(B,w)) = R(Y(Y(A,z-w)B,w)),
\label{eq:assoc}
\end{equation}
where the left-hand (resp.\ right-hand) side is the analytic
continuation from the domain $\left|z \right| > \left|w
\right|\/$ (resp.\ $ |w| > |z-w|\/$).
\end{prop}

This implies that for any $A, B \in V$ and $n \in \Z$, there exists $C_n
\in V$, such that $$\Res_{z=w} R(Y(A,z) Y(B,w)) (z-w)^n dz = Y(C_n,w).$$

\propref{assoc} also implies the following result.

\begin{prop}[\cite{FHL}, Prop.~3.5.1]    \label{converge}
For any $A_i \in V$, $i=1,\ldots,n$, the composition $Y(A_1,z_1) \ldots
Y(A_n,z_n)$ converges in the domain $|z_1| > \ldots > |z_n|$ and can be
continued to a rational operator-valued function on $(\C^\times)^3$ with
the poles only on the diagonals $z_i = z_j$.
\end{prop}

\section{Deformed Chiral Algebras}

\noindent {\bf Definition.} A {\cal deformed chiral algebra } (DCA) with a
representation is a collection of the following data:

\begin{itemize}
\item A vector space $V$ called the space of fields.

\item A vector space $W = \cup_{n\geq 0} W_n$ called the space of states,
which is union of finite-dimensional subspaces $W_n$. We consider a
topology on $W$ in which $\{ W_n \}_{n\geq 0}$ is the base of open
neighborhoods of $0$.

\item A linear map $Y: V\to \on{End} W \widehat{\otimes} [[z, z^{-1}]]$
that is for each $A \in V$, a formal power series $Y(A,z) = \sum_{n\in\Z}
A_n z^{-n}$, where each $A_n$ is a linear operator $W \arr W$, such that
$A_n \cdot W_m \subset W_{m+N(n)}, \forall m \geq 0$, for some $N(n) \in
\Z$.

\item A meromorphic function $S(x): \C^\times \arr \on{Aut}(V \otimes V)$,
satisfying the Yang-Baxter equation:
$$S_{12}(z) S_{13}(zw) S_{23}(w) =  S_{23}(w) S_{13}(zw) S_{12}(z),$$ for
all $z,w \in \C^\times$.

\item A lattice $L \subset \C^\times$, which contains the poles of 
$S(x)$.

\item A vector $\Omega \in V$, such that $Y(\Omega,z) = \on{Id}$.

\end{itemize}

The data of deformed chiral algebra should satisfy the following axioms:

\begin{itemize}

\item[(1)] For any $A_i \in V, i=1,\ldots,n$, the composition $Y(A_1,z_1)
\ldots Y(A_n,z_n)$ converges in the domain $|z_1| >> \ldots >> |z_n|$ and
can be continued to a meromorphic operator valued function $$R(Y(A_1,z_1)
\ldots Y(A_n,z_n)): (\C^\times)^n \arr \on{Hom}(W,\overline{W}),$$ where
$\overline{W}$ is the completion of $W$ with respect to its topology.

\item[(2)] Denote $R(Y(A,z)Y(B,w))$ by $Y(A \otimes B;z,w)$. Then
$$Y(A \otimes B;z,w) = Y(S(w/z)(B \otimes A);w,z).$$

\item[(3)] The poles of the meromorphic function $R(Y(A,z)Y(B,w))$ lie on
the lines $z = w\gamma$ where $\gamma\in L$. For each such line and $n\geq
0$, there exists $C_n \in V$, such that
$$\Res_{z=\gamma w} R(Y(A,z)Y(B,w)) (z-\gamma w)^n \dzz = Y(C_n,w).$$

\end{itemize}

\begin{rem} One can modify the last axiom by suppressing the condition
$n\geq 0$ and considering expansions at all $\gamma \in L$.

One can also allow more complicated topology on $W$. However, in all
examples that we give below, the topology of inductive limit of
finite-dimensional vector spaces on $W$ will suffice.\qed
\end{rem}

\begin{rem} According to \propref{converge}, an analogue of axiom (1) holds
in a VOA for all $n$, i.e., the case $n=2$ already implies the general $n$
case. We have not been able to prove that in a general DCA axiom (1) would
hold for all $n$ if it holds for $n=2$.\qed
\end{rem}

For $A, B \in V$, we write $$S(x)(A \otimes B) = \sum_i S_i(x)A \otimes
S^i(x)B.$$ Then axiom (2) can be expressed as $$R(Y(A,z)Y(B,w)) = \sum_i
R(Y(S_i(w/z)B,w)Y(S^i(w/z)A,z)).$$

Suppose we are given a DCA ${\cal V} = (V,W,Y,S(x),\Omega)$ and a subspace
$V' \subset V$, such that $V' \otimes V'$ is invariant under $S(x)$,
$\Omega \in V'$, and the following condition holds: for any $A,B \in V'$,
if $R(Y(A,z)Y(B,w))$ has a pole at $z=w\gamma$, then there are $C_n \in V',
n\geq 0$, such that $$\Res_{z=\gamma w} R(Y(A,z)Y(B,w)) (z-\gamma w)^n \dzz
= Y(C_n,w).$$ Then ${\cal V}' = (V',W,Y,S(x),\Omega)$ is a DCA, which we
will call a subalgebra of ${\cal V}$.

It is straightforward to define homomorphisms of deformed chiral algebras in
such a way that the kernel and the image are subalgebras.

The following result follows immediately from axioms (3) and (4) and the
residue formula, which should be compared with the Jacobi identity of
\cite{FLM,FHL}.

\begin{prop}
For $A, B, C \in V$ and any meromorphic function $f(z,w,u)$ on
$(\C^\times)^3$, and $a, b \in \C^\times$, the following identity holds:
\begin{align*}
&\Res_{z=ua} \Res_{w=ub} R(Y(A,z) Y(B,w) Y(C,u)) f(z,w,u) \dzz \dww \\ -
&\sum_i \Res_{w=ub} \Res_{z=ua} R(Y(S_i(w/z)B,w) Y(S^i(w/z)A,z) Y(C,u))
f(z,w,u) \dzz \dww \\ = &\Res_{w=ub} \Res_{z=wa/b} R(Y(A,z) Y(B,w) Y(C,u))
f(z,w,u) \dzz \dww.
\end{align*}
\end{prop}

\section{Examples of deformed chiral algebras}

Let $V = (V,Y,T,\Omega)$ be a super VOA and $W$ be a representation of
$V$. Then $(V,W,Y,(-1)^{p \otimes p},\Omega)$ where $p$ is the parity
operator on $V$, is a deformed chiral algebra.

The examples that we consider below are obtained via free field
realizations. We first construct a large DCA associated to an
infinite-dimensional Heisenberg algebra, and then define a chiral
subalgebra of this DCA, which is finitely generated in a certain sense.

\subsection{Deformed Heisenberg chiral algebra}

Fix a complex number $p$, such that $|p|<1$ and a complex number $q$.

Following \cite{A1,FF:w} consider the Heisenberg algebra ${\cal H}$ with
generators $\la_n, n \in \Z$, and relations
\begin{equation}    \label{hprel}
\left[ \la_n,\la_m \right] = -\frac{1}{n}
\frac{(1-q^n)(1-(p/q)^n)}{1+p^n} \delta_{n,-m}, \quad \quad n\neq 0,
\end{equation}
$$[\la_0,\la_n] = 0.$$

Let $\pi$ be the Fock representation of the algebra ${\cal H}$, which is
generated by a vector $v$, such that $\la_n v = 0, n\geq 0$. We define
a $\Z_+$--gradation on $\pi$ by setting $\deg \la_n = -n, \deg v =
0$.

Introduce the generating functions
$$\La^+(z) = :\La(z): = p^{-1/2} q^{- \la_0} :\exp \left( - \sum_{m\neq 0}
\la_m z^{-m} \right):,$$ $$\La^-(z) = :\La(zp^{-1})^{-1}: = p^{1/2}
q^{\la_0} :\exp \left( \sum_{m\neq 0} \la_m p^m z^{-m} \right):,$$ where
columns stand for the standard normal ordering. The Fourier coefficients of
$\La^\pm(z)$ are well-defined linear operators acting on $\pi$.

Denote $\theta_a(x) = \prod_{n>0} (1-xa^{n-1})(1-x^{-1}a^n)(1-a^n)$.
The following lemma is straightforward.

\vspace*{4mm}
\noindent {\bf Lemma-construction.} 
{\em Let $H$ be the vector space of polynomials in the variables
$\La^\pm_{i,n,m}, i,n,m \in \Z, i\geq 0$. Set $\Omega = 1$, and

$$Y(\La^{\ep_1}_{i_1,n_1,m_1} \ldots \La^{\ep_k}_{i_k,n_k,m_k},z)
= :\pa_z^{i_1} \La^{\ep_1}(z p^{n_1} q^{m_1}) \ldots \pa_z^{i_k}
\La^{\ep_k}(z p^{n_k} q^{m_k}):.$$

This operation gives rise to a structure of DCA on $(H,\pi)$ with $S(x): H
\otimes H \arr H \otimes H$ determined uniquely by axiom (2) and the Wick
rules. In particular,
$$S(\La^{\ep}_{0,n,m} \otimes \La^{\ep'}_{0,n',m'}) = S_{\La\La} \left(
\frac{w}{z} p^{n'-n} q^{m'-m} \right) \La^{\ep'}_{0,n',m'} \otimes
\La^{\ep}_{0,n,m},$$ where
\begin{equation}    \label{SLL}
S_{\La\La}(x) = \frac{\theta_{p^2}(xp)
\theta_{p^2}(xq^{-1}) \theta_{p^2}(xp^{-1}q)}{\theta_{p^2}(xp^{-1})
\theta_{p^2}(xq) \theta_{p^2}(xpq^{-1})}.
\end{equation}
}

\subsection{Deformed Virasoro chiral algebra}

Define the power series
$$T(z) = \sum_{m\in\Z} T_m z^{-m}$$ by the formula
\begin{equation}    \label{qvirrel}
T(z) = \La^+(z) + \La^-(z) = :\La(z): + :\La(zp^{-1})^{-1}:.
\end{equation}

The coefficients $T_n$ of the power series $T(z)$ are well-defined
linear operators on $\pi$. They satisfy the following
relations \cite{A1}:
\begin{equation}    \label{skao}
\sum_{l\geq 0} f_l \left(T_{n-l} T_{m+l} - T_{m-l} T_{n+l} \right) =
\frac{(1-q)(1-p/q)}{1-p} (p^{-n}-p^n) \delta_{n,-m},
\end{equation}
where $f_l$'s are given by the generating function
\begin{equation}    \label{fz}
f(x) = \sum_{l\geq 0} f_l x^l = \frac{1}{1-x}
\frac{(x|q,pq^{-1};p^2)_\infty}{(x|pq,p^2q^{-1};p^2)_\infty}.
\end{equation}

Here and below we use the following notation:
$$(a_1,a_2,\ldots,a_k;b)_\infty = \prod_{i=1}^k \prod_{n\geq 0} (1-a_i
b^n).$$

We associate to this algebra a subalgebra of the DCA $(H,\pi)$.

\vspace*{4mm}
\noindent {\bf Definition.} {\em The deformed chiral Virasoro algebra is the
smallest subalgebra $(\Vir,\pi)$ of the DCA $(H,\pi)$, which contains the
vector $T = \La^+ + \La^-$.}
\vspace*{4mm}

Thus, in particular, $Y(T,z) = T(z)$ is given by \eqref{qvirrel}.

\begin{prop}
There is a basis ${\cal B}$ in $\Vir$ containing $T$, with
respect to which $S(x)$ is diagonal, i.e.,
\begin{equation}    \label{exchange}
S(x)(A \otimes B) = S_{AB}(x) \cdot A \otimes B,
\end{equation}
where $S_{AB}(x)$ is a meromorphic function, for all $A,B \in {\cal
B}$. The functions $S_{AB}(x)$ are uniquely determined by the function
$S_{TT}(x) \equiv S_{\La\La}(x)$ given by \eqref{SLL}.
\end{prop}

\vspace*{4mm}
\noindent {\bf Proof.}
We can construct $\Vir$ inductively. For $|z| >> |w|$ we have
\begin{align}    \label{TTfusion}
T(z) T(w) &= f\wz^{-1} :\La(z) \La(w): + f\left( \frac{w}{zp} \right)
:\La(z) \La(wp^{-1})^{-1}: \\ \notag &+ f\left( \frac{wp}{z} \right)
:\La(zp^{-1})^{-1} \La(w): + f\wz^{-1} :\La(zp^{-1})^{-1}
\La(wp^{-1})^{-1}:,
\end{align}
where $f(x)$ is given by formula \eqref{fz}.

The analytic continuations satisfy \cite{FF:w}
\begin{equation}    \label{satisfy}
T(z) T(w) = \ST \left( \frac{w}{z} \right) T(w) T(z),
\end{equation}
where $\ST(x)$ is given by \eqref{SLL}. We will call a relation of this
kind an exchange relation. Note that the poles of the function $\ST(x)$ lie
in the lattice $L = p^{\Z} q^{\Z} \subset \C^\times$.

Next, according to axiom (3), $\Vir$ has to contain elements corresponding
to each pole in the right hand side of \eqref{TTfusion}. These are simple
poles occurring at the lines $z=wqp^2n, z=wq^{-1} p^{2n+1}, n\geq 0$, and
$z=wp^{\pm 1}$. The residue of $T(z) T(w)$ at the last two poles is a
constant, hence it corresponds to a multiple of $\Omega$. The residue at
each of the other poles can be computed explicitly (see examples below) and
it is a linear combination of normally ordered products of $\La^\pm(w
\gamma)$, where $\gamma \in L$. Let $R(w)$ be the residue, corresponding to
the line $z=wa$. We can compute the action of $S$ on the $T \otimes R$ and
$R \otimes R$, using the following formula.

\begin{lem}
Suppose we have three fields $A_i(z), i=1,2,3$, from a DCA $(V,W)$, which
satisfy $R(A_i(z) A_j(w)) = S_{ij}(w/z) R(A_j(w) A_i(z))$ for some
meromorphic functions $S_{ij}(x)$. Let $$B(w) = \Res_{u=w\gamma} R(A_2(u)
A_3(w)) \frac{du}{u}.$$ Then $A_1(z) B(w) = S_{12}(w\gamma/z) S_{13}(w/z)
B(w) A_1(z)$.
\end{lem}

This formula implies that if $A$ and $B$ belong to the set of basis
elements that we have already constructed, then $S_{AB}(x)$ is the product
of functions $\ST(x)$ with shifts by elements of $L$. Therefore the element
$R$ can be included in the basis ${\cal B}$.

We continue our construction inductively, by considering the OPEs of
the already constructed elements, and adding new elements to the basis
${\cal B}$ of $\Vir$ corresponding to the residues at various poles of
these OPEs. The fields corresponding to the new elements will be
linear combinations of normally ordered products of $\pa_z^i
\La^\pm(z)$ with shifts by elements of $L$, and they will have
exchange relations of the form \eqref{exchange} with the function
$S_{AB}(x)$ being the product of functions $\ST(x)$ with shifts by
elements of $L$.\qed

\subsection{Remarks}

\subsubsection{} The interesting question is to characterize the subalgebra
$\Vir$ inside the Heisenberg algebra $H$. In the undeformed case, this is
done by introducing the so-called screening operators. For general
$\W$--algebra $\W_\beta(\g)$, these are certain operators acting from the
VOA of $\W_\beta(\g)$ to its representations, corresponding to the simple
roots of $\g$, see \cite{FF:laws}. It is proved in \cite{FF:laws} that for
generic values of $\beta$ (i.e., generic central charge) the VOA of
$\W_\beta(\g)$ coincides with the intersection of kernels of the screening
operators.

The screening operators for the deformed $\W$--algebras associated to
$\sw_N$ have been constructed in \cite{FF:w,A2}; this construction is easy
to generalize to the case of an arbitrary simply-laced $\g$, see
\cite{FF:w}. The generators of the deformed $\W$--algebra $\W_{p,q}(\sw_N)$
commute with these operators. An important open question is to characterize
$\W_{p,q}(\sw_N)$ in terms of these operators.

\subsubsection{} The deformed $\W$--algebras have additive (or
``Yangian'') analogues, see \cite{Hou1}, which can be considered as
symmetries of the sine-Gordon and Toda theories. At the semi-classical
level, they correspond to the center of the extended Yangian double at the
critical level, see \cite{Hou2}. The formalism of DCA given in this section
also has an additive version, suitable for these algebras.

\subsubsection{}    \label{analogy}
The VOA of the $\W$--algebra can be extended by adding to it various
representations, i.e., various primary fields and their descendants. In
particular, in the case of the Virasoro algebra one can add the primary
fields $\Phi_{n,m}(z)$. It is well-known that the OPE of the fields
$\Phi_{1,2}^+(z)$ and $\Phi_{1,2}^-(z)$ contains the generating field
$T(z)$ of the Virasoro algebra.

It is interesting that the same pattern holds in the deformed case. Indeed,
the free field realization of the deformed fields $\Phi_{1,2}^\pm(z)$ has
already been given by Lukyanov and Pugay in \cite{LP1} (see also
\cite{Luk1}), and they play a prominent role in the subsequent work on the
RSOS models, see \cite{LP2}. Explicit computation shows that the OPE of the
deformed fields $\Phi_{1,2}^+$ and $\Phi_{1,2}^-$ contains the generating
field of the deformed Virasoro algebra. Furthermore, ``Yangian'' version of
this (scaling limit) can be interpreted as the appearance of the breather
particle in the scattering of soliton and anti-soliton of the sine-Gordon
theory, see \cite{M2,Luk2}.

Thus one obtains a remarkable connection between scattering in an
integrable quantum field theory and the OPE of a conformal field theory,
and correspondence $\Phi_{1,2}^{+(-)}(z) \leftrightarrow$ (anti-)soliton,
$T(z) \leftrightarrow$ breather, etc. Note however, that $z$ plays
completely different roles in the two theories: as momentum and coordinate
variables, respectively.

This example clearly indicates that the DCA of the deformed Virasoro
algebra (as well as other deformed $\W$--algebras) can, and perhaps should,
be extended by its ``primary fields''.

\subsubsection{} The notion of VOA has a semi-classical analogue:
``coisson algebra'' \cite{BD}, or ``vertex Poisson algebra'' \cite{EF}. It
is closely connected with the hamiltonian formalism of soliton theory. The
deformed chiral algebras also have a semi-classical analogue, which we plan
to discuss in a future paper.

\section{Deformed $\W$--algebra as a quantization of representation ring}

Let $\G$ be an affine Kac-Moody algebra. As we mentioned in the
introduction, a construction given in \cite{RS} associates to each
finite-dimensional representation $V$ of $U_p(\G)$ a generating series
$t_V(z)$ of central elements from $Z_p(\g)[[z,z^{-1}]]$.

This gives us a map ${\cal T}$ from the representation ring $\on{Rep}
U_p(\G)$ of $U_p(\G)$ to $Z_p(\g)$. Recall that there is a natural
$\C^\times$ action on $\on{Rep} U_p(\G)$: given a representation $\pi_V$ of
$U_p(\G)$ in a vector space $V$, and $z \in \C^\times$, we construct a new
representation $\pi_{V(z)})$ on the same vector space by the formula
$\pi_{V(z)}(a) = \pi_V(z^d a z^{-d})$, where $d$ is a generator of $\G$.

\begin{prop}[\cite{DE}] The map ${\cal T}: \on{Rep} U_p(\G) \arr
Z_p(\g)[[z,z^{-1}]]$ is a $\C^\times$--equivariant ring homomorphism,
i.e., $t_{U \oplus V}(z) = t_U(z) + t_V(z)$, $t_{U \otimes V}(z) = t_U(z)
t_V(z)$, and $t_{V(w)}(z) = t(zw)$.
\end{prop}

The construction of the map ${\cal T}$ is closely connected to the
construction of transfer-matrices of the integrable spin models
corresponding to $U_p(\G)$ (see \cite{RS}). More precisely, there is a
surjective homomorphism from the center $Z_p(\g)$ to the algebra of
transfer-matrices, under which $t_V(z)$ goes to the transfer-matrix
corresponding to the representation $V$. The transfer-matrices generate a
large commutative subalgebra of a non-commutative quantum algebra. The
discovery of \cite{FR} was that this algebra has a Poisson structure, which
can be further quantized (``secondary quantization''!). Thus,
transfer-matrices appear as quantum objects in one theory (they form a
commutative subalgebra in some non-commutative quantum algebra, e.g.,
Yangian or quantum affine algebra) and as classical objects in another
theory (the semi-classical limits of the generators of a non-commutative
symmetry algebra, or Faddeev-Zamolodchikov algebra). If the first theory is
rational, the second theory is trigonometric, and if the first is
trigonometric, then the second theory is elliptic.

Furthermore, in \cite{FR} we showed that the Baxter type formulas for the
eigenvalues of transfer-matrices computed by the Bethe Ansatz \cite{KS}
actually define a hamiltonian map from the algebra of transfer-matrices to
a Heisenberg-Poisson algebra (Miura transformation, or free field
realization). This hamiltonian map can also be quantized, at least in the
case of $\sw_N$ \cite{FF:w,A2}.

Now let us recall that in the undeformed case, the one-parameter family of
deformations of the center $Z(\g)$ is the $\W$--algebra $\W_\beta(\g)$
($\beta$ is the parameter related to the central charge, see
\cite{FF:laws}), and it actually has two classical limits: $\beta \arr 0$
and $\beta \arr \infty$. The first of them is $\W(\g)$ and the second is
$\W(^L\g)$, which is isomorphic to the center $Z(\g)$ of
$U(\G)_{\on{cr}}$. Thus, at ``one end'' of $\W_\beta(\g)$ one obtains
$Z(\g)$, and at ``the other end'' one obtains the Langlands dual $Z(^L\g)$.

We believe that the center $Z_p(\g)$ also has a one-parameter family of
deformations $Z_{p,q}(\g)$, which has two classical limits. One of them is
the limit $q \arr p$ in which $Z_{p,q}(\g)$ becomes $Z_p(\g)$ and the other
is $q \arr 1$.

Such a deformation has been constructed \cite{FF:w,A2} in the case when
$\g$ is $\sw_N$; in this case $^L\g$ is also $\sw_N$, and we know that in
the limit $q \arr 1$ we also obtain $Z_p(\sw_N)$. Note also that one
obtains $\W_\beta(\sw_N)$ in the limit $q \arr 1, p = q^\beta$ with $\beta$
fixed. The same pattern is probably true for all simply-laced $\g$ (see the
end of \cite{FF:w}).

The structure of the Poisson algebras $Z_p(\g)$ for arbitrary $\g$ and the
corresponding deformed Miura maps have been worked out in
\cite{FRK}. However, it is not clear how to quantize $Z_p(\g)$ when $\g$ is
not simply-laced. The naive expectation is that the $q \arr 1$ limit of the
quantum $\W$--algebra $Z_{p,q}(\g)$ is isomorphic to $Z_p(^L\g)$. However,
our explicit construction \cite{FRK} indicates that this is not so for non
simply laced $\g$. It is not clear to us at the moment what takes the role
of Langlands duality in this context.

If exists, $Z_{p,q}(\g)$ can be thought of as a quantization of $\on{Rep}
U_p(\G)$ (or the algebra of transfer-matrices in the corresponding spin
model). Let $V$ be a finite-dimensional representation of $U_p(\G)$ and
$t_V(z)$ be the corresponding generating series of elements from
$Z_p(\g)$. It is interesting to ask what is its deformation in
$Z_{p,q}(\g)$, which gives $t_V(z)$ in the limit $q \arr p$ (of course,
such a deformation is not unique, since we can always add $(1-p/q)$ times
something).

In the case when $\g=\sw_N$, the only case, in which the quantum
$\W$--algebra hs been worked out explicitly, we only know what are the
deformations of $t_V(z)$, where $V$ is a fundamental representation. What
we would like to do is to express deformation of other $t_V(z)$ in terms of
these. That's where the concept of DCA becomes useful.

In the next section we will give examples of such expressions in the case
when $\g=\sw_2$. In this case, for each positive integer $n$ there exists
an irreducible $(n+1)$--dimensional representation. Denote the
corresponding series by $t_n(z)$. This series can be written in terms of
the deformed Miura transformation as follows:
\begin{align*}
t_n(z) &= \La^+(w) \La^+(wp) \La^+(wp^2) \ldots \La^+(wp^n) \\
&+ \La^-(w) \La^+(wp) \La^+(wp^2) \ldots \La^+(wp^n) \\
&+ \La^-(w) \La^-(wp) \La^+(wp^2) \ldots \La^+(wp^n) \\
&+ \ldots \\
&+ \La^-(w) \La^-(wp) \La^-(wp^2) \ldots \La^-(wp^n),
\end{align*}

The deformation of $t_1(z)$ in $\W_{p,q}(\sw_2)$ is $T_1(z) \equiv
T(z)$. It turns out that the deformation $T^q_n(z)$ of the series $t_n(z)$
can be obtained inductively as a residue occurring under a multiple fusion
of the field $T(z)=T_1(z)$.

\section{Examples}

\subsection{Examples of operator product expansions}

We obtain from formula \eqref{TTfusion}:
$$\Res_{z=wq} T(z) T(w) \dzz = \frac{(1-q^{-1})(p,p^2
q^{-1};p^2)_\infty}{(p^2,pq^{-1};p^2)_\infty} \times$$ $$\left( :\La(w)
\La(wq): + :\La(wp^{-1}q)^{-1} \La(wp^{-1})^{-1}: +
\frac{(1-p)(1+q)}{(1-pq)} :\La(wq) \La(wp^{-1})^{-1}: \right),$$

Denote by $T_2^{q}(w)$ the field appearing in brackets in the right hand
side. Its limit when $q \arr p$ equals $$t_2(w) = \La(w) \La(w) +
\La(w)^{-1} \La(wp^{-1})^{-1} + \La(wp) \La(wp^{-1})^{-1}.$$

We find by induction that
\begin{align*}
\Res_{z=wq^n} & T(z) T_n^q(w) \dzz = \prod_{j=1}^n
f(q^{-n+j-1})^{-1} \times
\\ &\left( c^0_{n+1} :\La^+(w) \La^+(wq) \La^+(wq^2) \ldots \La^+(wq^n):
\right.
\\ &+ c^1_{n+1} :\La^-(w) \La^+(wq) \La^+(wq^2) \ldots \La^+(wq^n):
\\ &+ c^2_{n+1} :\La^-(w) \La^-(wq) \La^+(wq^2) \ldots \La^+(wq^n):
\\ &+ \ldots
\\ &\left. + c^{n+1}_{n+1} :\La^-(w) \La^-(wq) \La^-(wq^2) \ldots
\La^-(wq^n): \right),
\end{align*}
where $f(x)$ is given by \eqref{fz}, and the numbers $c^i_n$ can be
determined inductively by the formulas $c^0_n=c^n_n = 1$ and $c^i_{n+1} =
c^i_n \prod_{j=1}^i g(q^{n-j+1})$ for $0<i<n+1$, where $$g(x) =
\frac{(1-xq)(1-xp/q)}{(1-x)(1-xp)}.$$ We denote the term in brackets by
$T^q_{n+1}(w)$. From the point of view of the analogy discussed in
\remref{analogy}, $T_n^q(z)$ corresponds to the $n$--breather particle of
the sine-Gordon model.

The limit of $T^q_n(z)$ when $q \arr p$ equals $t_n(z)$. Thus, the field
$T^q_n(z)$ can be viewed as a deformation of the transfer-matrix
corresponding to the $(n+1)$--dimensional representation of
$U_p(\widehat{\sw}_2)$. It is interesting if one can see the tensor product
structure of representations of $U_p(\widehat{\sw}_2)$ in the scattering of
$n$--breathers.

Similarly, we can construct the deformations $T_n^{p/q}(z)$ by considering
the residues at $z=w(p/q)^{n-1}$. In particular,
\begin{align*}
T_2^{p/q}(w) &= :\La(w) \La(wpq^{-1}): \\ & + :\La(wp^{-1})^{-1}
\La(wq^{-1})^{-1}: + \frac{(1-p)(1+pq^{-1})}{1-p^2q^{-1}} :\La(wpq^{-1})
\La(wp^{-1})^{-1}:.
\end{align*}
The field $T^{p/q}_n(z)$ becomes $t_n(z)$ in the limit $q \arr 1$.

\subsection{Formal power series relations}

Recall that one can use the OPE in the vertex operator algebra to
compute the commutation relations between the Fourier coefficients of
the fields. These relations can be written in the form of power series
relations. Explicitly, given two fields $A(z)$ and $B(w)$ in a VOA
with the singular part of their OPE of the form $$A(z) B(w) =
\sum_{n>0} \frac{C_n(w)}{(z-w)^n},$$ we can derive the following
relation on the formal power series (see, e.g., \cite{Kac}):
$$[A(z),B(w)] = \sum_{n>0} \frac{1}{(n-1)!} C_n(w) (w \pa_w)^{n-1}
\delta \left( \frac{w}{z} \right).$$ In deriving this formula, one uses
the locality $R(A(z) B(w)) = R(B(w) A(z))$ in an essential way.

The locality axiom becomes more complicated in the case of
DCA. Nevertheless, it is possible to derive certain quadratic-linear
relations on the Fourier coefficients of the fields as well. Suppose we
have two fields $A(z)$ and $B(w)$, which satisfy $$R(A(z) B(w)) = S \left(
\frac{w}{z} \right) R(B(w) A(z)).$$ Let us factorize the function $S(x)$
into a product of two functions, analytic in the interior and the exterior
of a circle $|x| = \al$ (Riemann problem): $$S(x) = g_+(x)^{-1} g_-(x).$$
Here $g_+(x)$ can be expanded into a power series in $x$, while $g_-(x)$
can be expanded into a power series in $x^{-1}$. It is then easy to show,
in the same way as in the case of VOA, that the difference $$g_+ \left(
\frac{w}{z} \right) A(z) B(w) - g_- \left( \frac{z}{w} \right) B(w) A(z),$$
considered as a formal power series in $z$ and $w$, is the sum of {\em
finitely many} terms of the form $$C(w) (w \pa_w)^k \delta\left(
\frac{w}{z} \gamma \right),$$ where
$$C(w) = \on{Res}_{z=w\gamma} g_+ \wz (z-w)^k R(A(z) B(w)) \frac{dz}{z}$$
up to a scalar. This way one obtains quadratic-linear relations on the
Fourier coefficients of fields in a DCA (compare with Sect.~7 of
\cite{FF:w}). But note that if we solve the Riemann problem differently, or
on a different circle, we may obtain a different relation. This is
illustrated by the examples below.

\begin{rem} We want to use this opportunity to correct a typo in formula
(6.2) of \cite{FF:w}: $\displaystyle f_{m,N} \left( \frac{z}{w} \right)$
there should be replaced by $\displaystyle f_{m,N} \left(
\frac{z}{wp^{m-1}} \right)$.
\end{rem}

\subsubsection{Relations between $T_1(z)$ and $T_1(w)$}
We have already seen relations between the Fourier coefficients of $T_1(z)$
and $T_1(w)$ in formula \eqref{skao}. Those relations can be written as
relations on power series:
\begin{multline}    \label{odin}
f \left(\frac{w}{z}\right) T(z) T(w) - f \left(\frac{z}{w}\right)
T(w) T(z) \\ = \frac{(1-q)(1-p/q)}{1-p} \left( \delta \left( \dfrac{w}{zp}
\right) - \delta \left( \dfrac{wp}{z} \right) \right),
\end{multline}
where $f(x)$ is given by formula \eqref{fz}. From the point of view of the
general formalism outlined above, this relation corresponds to the
factorization $$\ST(x) = f(x)^{-1} f(x^{-1})$$ (cf. \cite{FF:w}). But there
are many other factorizations, for instance,
$$\ST(x) = f(xp) f(x^{-1} p^{-1})^{-1}.$$ It results in a new relation:
$$f \left( \frac{wp}{z} \right)^{-1} T_1(z) T_1(w) - f \left( \frac{z}{wp}
\right)^{-1} T_1(w) T_1(z)$$ $$= \frac{(1-q^{-1})(1-pq^{-1})}{1-pq^{-2}}
\left( T_w^q(w) \delta \left( \frac{wq}{z} \right) - T_2^{p/q}(w) \delta
\left( \frac{wp}{zq} \right) \right)$$ $$+
\frac{(1-q)(1-pq^{-1})(1-p^2)}{(1-pq)(1-p^2q^{-1})} \delta \left(
\frac{w}{zp} \right).$$

\subsubsection{Relations between $T_1(z)$ and $T^q_2(w)$}
We have a factorization: $$S_{T_1,T^q_2}(x) = \ST(x) \ST(xq) = F_1(x)^{-1}
F_1(x^{-1} q^{-1}),$$ where $$F_1(x) =
\frac{(xp,xp^2q,xpq^{-1},xq^2;p^2)_\infty}{(x,xpq,xp^2q^{-1},xpq^2;
p^2)_\infty}.$$ The corresponding relation is
$$F_1 \wz T_1(z) T_2^q(w) - F_1 \left( \frac{z}{wq} \right) T^q_2(w) T_1(z)
=$$ $$\frac{(1-pq^{-1})(1 - q^2)}{1-pq} \left( T_1(wq) \delta \left(
\frac{w}{zp} \right) - T_1(w) \delta \left( \frac{wpq}{z} \right)
\right).$$

\subsubsection{Relations between $T^q_2(z)$ and $T^q_2(w)$}
$$F_2 \wz T_2^q(z) T_2^q(w) - F_2 \left( \frac{z}{w} \right) T_2^q(w)
T_2^q(z) =$$ $$\frac{(1-p)(1+q)(1-q^2)(1-pq^{-1})}{(1-pq)^2} \left(
\wt{T}(wq) \delta \left( \frac{w}{zpq} \right) - \wt{T}(w) \delta \left(
\frac{wpq}{z} \right) \right)$$ $$+
\frac{(1-q)(1-pq^{-1})(1-pq^{-2})(1-q^2)}{(1-p)(1-q^{-1})(1-pq)} \left(
\delta \left( \frac{w}{zp} \right) - \delta \left( \frac{wp}{z} \right)
\right),$$ where $$F_2(x) = \frac{1}{1-x}
\frac{(xp^2q,xpq^{-1},xpq^{-1},xq,xq^2,xpq^{-2};p^2)_\infty}{(xpq,xpq,
xq^{-1},xpq^2,xp^2q^{-2},xp^2q^{-1};p^2)_\infty},$$
and
\begin{align*}
\wt{T}(w) &= \Res_{z=wpq^2} f \wz \frac{T(z) T(w)}{z-wpq^2} \; dz \\ &=
:\La(w) \La(wpq^2): + \frac{(1-pq)(1-q^3)}{(1-pq^2)(1-q^2)} :\La(w)
\La(wq^2)^{-1}: \\ &+ \frac{(1-pq^3)(1-p^2q)}{(1-pq^2)(1-p^2q^2)}
:\La(wp^{-1})^{-1} \La(wpq^2): + :\La(wp^{-1})^{-1} \La(wq^2)^{-1}:.
\end{align*}

The Fourier coefficients of $\wt{T}(w)$ can be expressed quadratically in
terms of the Fourier coefficients of $T(w)$. Indeed, let $C$ be a circle of
large radius on the $z$ plane, and $C'$ be a circle of small radius. Then 
\begin{equation}    \label{manip}
\int_{C} f \wz \frac{R(T(z) T(w))}{z-wpq^2} \; dz - \int_{C'} f \left(
\frac{w}{z} \right) \frac{R(T(z) T(w))}{z-wpq^2} \; dz
\end{equation}
(here $f(x)$ is considered as a meromorphic function and not as a power
series) equals the sum of residues of the integrand. But the integrand
equals
\begin{align*}
\frac{1}{z-wpq^2} \left( :\La(z) \La(w): + \frac{(z-w q^{-1})(z-w p^{-1}
q)}{(z-w)(z-w p^{-1})} :\La(z) \La(wp^{-1})^{-1}: \right. \\ \notag
\left. + \frac{(z-w q)(z-w pq^{-1})}{(z-w)(z-w p)} :\La(zp^{-1})^{-1}
\La(w): + :\La(zp^{-1})^{-1} \La(wp^{-1})^{-1}: \right).
\end{align*}
Hence the poles occur only when $z=wpq^2$, $z=wp^{\pm1}$, and $z=w$. The
contribution of the first coincides with $\wt{T}(w)$, while the remaining
poles give us a constant $$-q^2
\frac{(1-q)(1+p)(1-pq^{-1})}{(1-q^2)(1-p^2q^2)}.$$

On the other hand, since $$f \wz R(T(z) T(w)) = f \left( \frac{z}{w}
\right) R(T(w) T(z))$$ by formula \eqref{satisfy}, we can rewrite
\eqref{manip} as
\begin{equation}    \label{manip1}
\int_{C} f \wz \frac{R(T(z) T(w))}{z-wpq^2} \; dz - \int_{C'} f \left(
\frac{z}{w} \right) \frac{R(T(w) T(z))}{z-wpq^2} \; dz.
\end{equation}
But now in the first summand we have $|z| >> |w|$, so we can strip
$R(\cdot)$ and consider the integrand as a power series in $w/z$, while in
the second summand we have $|w| >> |z|$, and we can consider it as a power
series in $z/w$. Hence \eqref{manip1} equals
$$\sum_{k \in \Z} \left( \sum_{i\leq 0} T_i T_{k-i} \al_i +
\sum_{i>0} T_{k-i} T_i \al_i \right) w^{-k},$$ where $$\sum_{i\leq 0} \al_i
x^{-i} = \frac{f(x)}{1-xpq^2},$$
$$\sum_{i>0} \al_i x^{-i} = \frac{f(x^{-1})}{1-xpq^2}.$$

Thus, we obtain finally:
$$\wt{T}(w) = \sum_{k \in \Z} \left( \sum_{i\leq 0} T_i T_{k-i} \al_i +
\sum_{i>0} T_{k-i} T_i \al_i \right) w^{-k} + q^2
\frac{(1-q)(1+p)(1-pq^{-1})}{(1-q^2)(1-p^2q^2)}.$$

\end{document}